\documentclass[reqno,12pt]{amsart}
\usepackage{a4wide,color,eucal,enumerate,mathrsfs}
\usepackage[normalem]{ulem}
\usepackage{psfrag}          
\usepackage{amsmath,amssymb,epsfig}
\numberwithin{equation}{section}

\usepackage[latin1]{inputenc}

\usepackage[pdfborder={0 0 0}]{hyperref}

\newtheorem{theorem}{Theorem}[section]

\newtheorem{question}[theorem]{Question}

\newtheorem{definition}[theorem]{Definition}

\newtheorem{remark}[theorem]{Remark}

\newcommand{\CD}{{\sf CD}}
\newcommand{\RCD}{{\sf RCD}}
\newcommand{\MCP}{{\sf MCP}}

\newcommand{\ppi}{{\mbox{\boldmath$\pi$}}}

\newcommand{\diam}{\mathop{\rm diam}\nolimits}

\renewcommand{\d}{{\mathrm d}}

\newcommand{\R}{\mathbb{R}}

\newcommand{\mm}{\mathfrak m}
\newcommand{\nn}{\mathfrak n}
\newcommand{\sfd}{{\sf d}}

\newcommand{\geo}{{\rm{Geo}}}

\newcommand{\restr}[1]{\lower3pt\hbox{$|_{#1}$}}

\newcommand{\p}{{\mathtt p}}

\begin{document}

\title[Failure of topological rigidity results for $MCP(K,N)$]
{Failure of topological rigidity results for 
the measure contraction property}

\author{Christian Ketterer}
\address{Institut f\"ur Angewandte Mathematik, Universit\"at Bonn\\
Endenicher Allee 60 \\
D-53115 Bonn \\
Germany 
}
\email{ketterer@iam.uni-bonn.de}

\author{Tapio Rajala}
\address{Department of Mathematics and Statistics\\
P.O.Box 35 (MaD), FI-40014 University of Jyv\"askyl\"a\\ Finland
}
\email{tapio.m.rajala@jyu.fi}


\begin{abstract}
 We give two examples of metric measure spaces satisfying the measure contraction property $\MCP(K,N)$
 but having different topological dimensions at different regions of the space.
 The first one satisfies $\MCP(0,3)$ and contains a subset isometric to $\R$, but does not topologically split.
 The second space satisfies $\MCP(2,3)$ and has diameter $\pi$, which is the maximal possible diameter for 
 a space satisfying $\MCP(N-1,N)$, but is not a topological spherical suspension. The latter example
 gives an answer to a question by Ohta.
\end{abstract}

\maketitle


\section{Introduction}

There are several generalizations of Ricci curvature lower bounds to metric measure spaces that use optimal mass transportation.
The definitions that have received the most attention are the $\CD(K,N)$ definitions by Sturm \cite{S2006I, S2006II}, and Lott and Villani
\cite{LV2009}. In the definitions the first number $K \in \R$ always refers to a lower bound on the Ricci curvature
and the second number $N \in [1,\infty]$ refers to an upper bound on the dimension of the space. 
A less restrictive version of Ricci curvature lower bounds, called the measure contraction property $\MCP(K,N)$
was studied by Sturm \cite{S2006II}, and by Ohta \cite{o2007a, o2007b}. Recently also a more restrictive class of definitions
called the Riemannian Ricci curvature lower bounds $\RCD^*(K,N)$ have been considered \cite{AGS2011, AGMR2012, Erbar-Kuwada-Sturm13, AMS2013}.

Common to all these definitions is stability under the measured Gromov-Hausdorff convergence. However, with the $\CD(K,N)$ and
$\MCP(K,N)$ definitions this stability has sometimes been sacrificed by assuming the space to be nonbranching. For the $\RCD^*(K,N)$ definitions
such extra assumption seems 
 unnecessary. This can be seen from the fact that such spaces are already essentially nonbranching \cite{RS2012, GRS2013}.
For some time it was not clear if the nonbranching assumption was just a technical assumption that was made to simplify the proofs.
This was the case for some of the results, like the local Poincar\'e inequalities \cite{R2012a, R2011}.
However, it was recently shown by the second named author that in order to obtain a local-to-global property for $\CD(K,N)$ 
it is necessary to assume nonbranching \cite{R2013b}. Notice that, as our above discussion suggests, for the $\RCD^*(K,N)$ 
a nonbranching assumption is not needed for the local-to-global property \cite{AMS2013b, GRS2013}.

In this paper we further study the necessity of the nonbranching assumption in the theory of
Ricci curvature lower bounds in metric measure spaces. This time we focus on geometric
(or topological) rigidity results. More specifically, we look at two classical rigidity results
from Riemannian geometry.
The first result is the Cheeger-Gromoll splitting theorem \cite{CG1971} saying
that a Riemannian manifold $\mathbf{M}$ with non-negative Ricci curvature containing an infinite line-segment
is isometric to $\R \times \mathbf{N}$ where $\mathbf{N}$ is again a Riemannian manifold with non-negative Ricci curvature.
The second rigidity result is the Cheng maximal diameter theorem \cite{C1975}
stating that an $n$-dimensional Riemannian manifold with Ricci curvature bounded from below by $n-1$
and with diameter $\pi$ is necessarily the standard sphere $\mathbb{S}^n$. Recall that,
by the Bonnet-Meyers theorem, $\pi$ is the maximal possible diameter for an $n$-dimensional Riemannian 
manifold with Ricci curvature bounded from below by $n-1$.

Both of the above mentioned rigidity results have been proven in the abstract setting of $\RCD^*(K,N)$ spaces.
The Cheeger-Gromoll splitting theorem was generalized to $\RCD(0,N)$ spaces by Gigli \cite{G2013}.
Using this generalization and by studying the metric cones over $\RCD^*(K,N)$ spaces the first named author recently generalized
the Cheng maximal diameter theorem to $\RCD^*(K,N)$ spaces \cite{K2013b}.

Now that we know that the rigidity results hold for $\RCD^*(K,N)$ spaces it is natural to ask if they hold
for the less restrictive definitions $\CD(K,N)$ and $\MCP(K,N)$.
Since $\CD(K,N)$ spaces include nonRiemannian Finsler manifolds, it would be unreasonable to
expect the Cheeger-Gromoll and Cheng theorems to hold in the same sharp form as for the $\RCD^*(K,N)$ spaces.
However, Ohta has shown that topological versions of the rigidity results still hold in the Finsler setting.
Indeed, he has proven a diffeomorphic splitting theorem for Finsler manifolds of
nonnegative weighted Ricci curvature \cite{O2013b}.
Ohta also proved in \cite{o2007b} that in nonbranching $\MCP(K,N)$ spaces 
a maximal diameter theorem holds in the following topological form. 

\begin{theorem}[Ohta \cite{o2007b}]\label{thm:ohta}
 Let $(X,\sfd,\mu)$ be a compact metric measure space satisfying the $\MCP(N-1,N)$ property for some $N >1$
 and assume that there exist $x_N, x_S \in X$ with $\sfd(x_N,x_S) = \pi$ such that
 \begin{equation}\label{eq:cutlocus}
   \textrm{Cut}(x_N) \setminus \{x_S\} = \textrm{Cut}(x_S) \setminus \{x_N\} = \emptyset.
 \end{equation}
 Then there exists a topological measure space $(Y,\nu)$ such that $(X,\mu)$ is
 the spherical suspension of $(Y,\nu)$ as a topological measure space.
\end{theorem}

Here the \emph{cut locus} $\textrm{Cut}(x)$ of a point $x \in X$ is the set of points $z \in X$ such that there
exist at least two distinct minimal geodesic between $x$ and $z$. Under the nonbranching assumption
\eqref{eq:cutlocus} is satisfied.
In \cite{o2007b, O2013} Ohta asked if the topological maximal diameter result holds for $\MCP(K,N)$ spaces
without the nonbranching assumption. We show that this is not the case.

\begin{theorem}\label{thm:2}
 There exists a compact geodesic metric measure space $(Y,\sfd,\nn)$ satisfying $\MCP(2,3)$ and $\diam(Y) = \pi$.
 Still the space is not a topological spherical suspension.
\end{theorem}

A slightly easier construction than the one showing the failure of the topological maximal diameter theorem
shows that also topological splitting results fail for $\MCP(0,N)$ spaces.

\begin{theorem}\label{thm:1}
 There exists a complete geodesic metric measure space $(X,\sfd,\mm)$ satisfying $\MCP(0,3)$ and containing an
 isometrically embedded copy of $\R$. Still the space does not topologically split.
\end{theorem}

Still several questions remain open
between the positive results and the counter-examples stated in Theorems \ref{thm:2} and \ref{thm:1}. Two obvious questions are the following.

\begin{question}\label{q:1}
 Does a topological maximal diameter theorem hold in $\CD(K,N)$ spaces?
\end{question}

\begin{question}
 Does a topological splitting theorem hold in $\CD(0,N)$ spaces?
\end{question}

Notice that there exist spaces satisfying $\CD(K,N)$ with positive $K$ that contain lots of branching geodesics \cite{O2013},
so Question \ref{q:1} does not have an obvious answer via Theorem \ref{thm:ohta}.

These questions are also related to the local structure of spaces with different Ricci curvature lower bounds.
It is still unknown even for $\RCD^*(K,N)$ spaces if the Hausdorff dimension or any other relevant dimension
is constant almost everywhere. Perhaps the currently best positive result in this direction is that $\RCD^*(K,N)$ spaces
have euclidean weak tangents at almost every point \cite{GMR2013}. The examples of this paper show that the least 
restrictive definitions of Ricci curvature lower bounds discussed here - the $\MCP(K,N)$ conditions - do not imply that the Hausdorff (or topological)
dimension of the space is almost everywhere constant.

\section{Preliminaries}

Before defining the spaces mentioned in Theorems \ref{thm:2} and \ref{thm:1} we recall some basics
including the definition of $\MCP(K,N)$ and its connection to other definitions of Ricci curvature lower bounds in metric measure spaces.

First of all, by a geodesic $\gamma$ in a metric space $(X,\sfd)$ we mean a map $\gamma \colon [0,1] \to X$
satisfying 
\[
\sfd(\gamma(t),\gamma(s)) = |s-t|\sfd(\gamma(0),\gamma(1)) \qquad \text{ for all }s,t \in [0,1].            
\]
The set of
all geodesics in $(X,\sfd)$ is denoted by $\geo(X)$.
For $\gamma \in \geo(X)$ we abbreviate $\gamma_t := \gamma(t)$ for all $t \in [0,1]$.
We also write length as $l(\gamma) := \sfd(\gamma_0,\gamma_1)$ for all $\gamma \in \geo(X)$
and define the evaluation maps
\[
e_t \colon \geo(X) \to X \colon \gamma \mapsto \gamma_t \qquad\text{ for all }t \in [0,1]. 
\]

The measure contraction property $\MCP(K,N)$ is defined using the functions
$\mathbf{s}_K$ with $K \in \R$. They are defined as
\begin{align*}
 \mathbf{s}_K \colon [0,\pi/\sqrt{K}) \to \R \colon t \mapsto \frac{\sin(\sqrt{K}t)}{\sqrt{K}},& \qquad \text{ for }K> 0,\\
 \mathbf{s}_K \colon [0,\infty) \to \R \colon t \mapsto \frac{\sinh(\sqrt{-K}t)}{\sqrt{-K}},& \qquad \text{ for }K< 0 \text{ and}\\
 \mathbf{s}_K \colon [0,\infty) \to \R \colon t \mapsto t,& \qquad \text{ for }K= 0. 
\end{align*}

Using these we define for all $t \in [0,1]$ and $d \ge 0$
\[
 \varsigma_{K,N}^{(t)}(d) = t\left(\frac{\mathbf{s}_K(td/\sqrt{N-1})}{\mathbf{s}_K(d/\sqrt{N-1})}\right)^{N-1}
\]
if $N > 1$ and $K \in \R$, and $\varsigma_{K,1}^{(t)}(d) = t$ if $K \le 0$.

Now we can introduce the $\MCP(K,N)$ condition as defined by Ohta \cite{o2007a}.

\begin{definition}
 For $K \in \R$, $N>1$, or $K \le 0$ and $N=1$ a metric measure space $(X,\sfd,\mm)$ is
 said to satisfy the \emph{$(K,N)$-measure contraction property}, $\MCP(K,N)$ for short,
 if for every point $x \in X$ and a measurable set $A \subset X$ with $0 < \mm(A) < \infty$
 (and $A \subset B(x,\pi\sqrt{(N-1)/K})$ is $K>0$) there exists a probability measure $\ppi$ on $\geo(X)$
 such that
 \begin{enumerate}
  \item We have $(e_0)_\sharp\ppi = \delta_x$ and $(e_1)_\sharp\ppi = \mm(A)^{-1}\mm\restr{A}$;
  \item For $t \in [0,1]$,
        \[
         \d\mm \ge (e_t)_\sharp\left(\varsigma_{K,N}^{(t)}(l(\gamma))\mm(A)\d\ppi(\gamma)\right)
        \]
        holds as measures on $X$.
 \end{enumerate}
\end{definition}

Let us then briefly give the connection of the $\MCP(K,N)$ condition to other similar curvature-dimension bounds.
The $\CD(K,N)$ conditions by Sturm \cite{S2006I, S2006II} and Lott and Villani \cite{LV2009} use the same
weights $\varsigma_{K,N}^{(t)}$ as $\MCP(K,N)$ but require a convexity type inequality between any two probability
measures along some optimal transport geodesic. 
As shown by the second named author in \cite{R2011}, $\CD(K,N)$ implies $\MCP(K,N)$.
There is also a more strict version of $\MCP(K,N)$ by Sturm \cite{S2006II}. He required the choice of 
geodesics be given by a Markov kernel. It is not known if the $\CD(K,N)$ conditions imply this more strict version
of $\MCP(K,N)$. In this paper we will work with the definition of Ohta.

Changing slightly the weights $\varsigma_{K,N}^{(t)}$ to a more PDE-friendly versions give rise
to so-called reduced definitions that usually carry a '$*$' in their abbreviation. The first set of
reduced definitions, $\CD^*(K,N)$ were considered by Bacher and Sturm in \cite{BS2010}. In
nonbranching spaces it was shown by Cavalletti and Sturm \cite{Cavalletti-Sturm12} that $\CD^*(K,N)$ actually implies $\MCP(K,N)$.
Recently Cavalletti has shown \cite{Cavalletti12} that in nonbranching spaces $\CD^*(K,N)$ gives the required convexity-type inequalities
for the $\CD(K,N)$ condition for a much larger class of transports. It is still open if the $\CD^*(K,N)$ really 
self-improves to $\CD(K,N)$.

The most recent additions to the set of definitions are the Riemannian Ricci curvature lower bounds \cite{AGS2011, AGMR2012, Erbar-Kuwada-Sturm13, AMS2013}.
These definitions start from the reduced weights and are therefore denoted by $\RCD^*(K,N)$. As the
word 'Riemannian' in the name of the condition suggests, the $\RCD^*(K,N)$ condition rules out nonRiemannian
Finsler structures that the $\CD^*(K,N)$ condition includes. As was noted in \cite{GMR2013} the $\RCD^*(K,N)$ conditions
imply $\MCP(K,N)$ without the nonbranching assumption that was assumed in \cite{Cavalletti-Sturm12}.
As in the case of $\CD(K,N)$, it is not known if $\RCD^*(K,N)$ self-improves to $\RCD(K,N)$.

There are also other curvature-dimension conditions besides the ones mentioned so far. One interesting notion is the
coarse Ricci curvature by Ollivier \cite{Ollivier2009}. Because by \cite{o2007a} we know that $\MCP(K,N)$ spaces
satisfy a Bishop-Gromov volume comparison theorem, the $\MCP(K,N)$ spaces (and in particular our examples) have
a lower bound on the coarse Ricci curvature associated to the $r$-step random walk, see \cite{Kitabeppu2013}.

\section{Examples of $\MCP(K,N)$ spaces with non-constant dimension}

We now turn to the examples stated in Theorems \ref{thm:2} and \ref{thm:1}.
They are both constructed as closed connected subsets of $\R^2$ equipped with the distance $\sfd$ coming from the $l^\infty$-norm,
i.e.
\begin{equation}\label{eq:dist}
 \sfd((x_1,y_1),(x_2,y_2)) = \sup(|x_1-x_2|,|y_1-y_2|) \qquad \text{for all }(x_1,y_1),(x_2,y_2) \in \R^2. 
\end{equation}

The first example we present here is a space satisfying $\MCP(0,3)$.
It is a tangent space of the second space satisfying $\MCP(2,3)$. Therefore, instead of
proving the $\MCP(0,3)$ property, we could refer to this fact. (See Section \ref{sec:remarks} for more discussion.)
We will still prove the $\MCP(0,3)$ property in the first example because it is simpler and
makes it easier to understand the second example.

In checking the examples we will use the convention that we are transporting from a dirac mass at
$(\tilde x,\tilde y)$ to a measure uniformly distributed on a set $A$.

\subsection{Failure of topological splitting}
\begin{figure}
  \centering
   \includegraphics[width=0.4\textwidth]{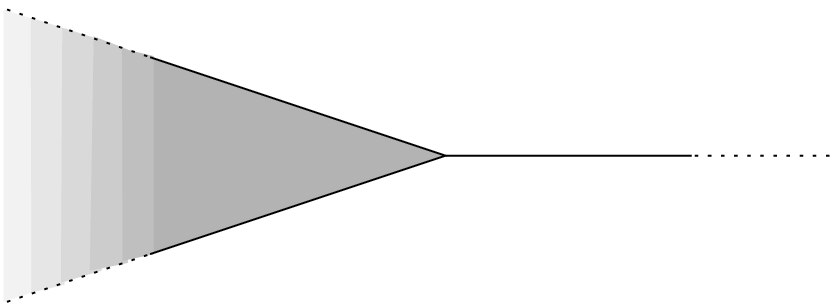}
  \caption{The space $(X,\sfd,\mm)$ satisfying $\MCP(0,3)$.}
  \label{fig:example1}
 \end{figure}

Let us start with the easier example of Theorem \ref{thm:1} of a space $(X,\sfd,\mm)$ satisfying $\MCP(0,3)$ containing an isometric embedding
of the euclidean real-line. As a subset of $\R^2$ the space $X$ is defined as
\[
 \{(x,y) \in \R^2 \,:\,  x \le -3|y| \le 0\}  \cup (\R_+ \times \{0\}),
\]
see Figure \ref{fig:example1} for an illustration.
 
The distance $\sfd$ on the space is the $l^\infty$-distance \eqref{eq:dist}.
The measure $\mm$ restricted to the half-infinite line segment $L = (\R_+ \times \{0\})$ is $\mathcal{H}^1\restr{L}$,
where $\mathcal{H}^1$ is the one-dimensional Hausdorff measure. On the cone $C = \{(x,y) \in \R^2 \,:\, x \le -3|y| \le 0\}$
we define $\mm$ as a weighted Lebesgue measure via
\[
 \d\mm = \frac{3}{-2x}\d\mathcal{L}_2.
\]
In other words, the projection of $\mm$ to the first coordinate is the Lebesgue measure $\mathcal{L}_1$
and the measure $\mm$ is equally distributed in $X \cap (\{x\}\times\R)$ for all $x \in \R$.

 
Let us then check that the space $(X,\sfd,\mm)$ satisfies $\MCP(0,3)$. We will always transport measure linearly in the horizontal
direction. Hence it is sufficient to check $\MCP(0,3)$ separately for sets $A_{\mathbf{x}} = \{(\mathbf{x},y) \in A \,:\, y \in \R\}$ with $\mathbf{x} \in \R$.
There are several cases to check depending on the position of the measure $\mu_0 = \delta_{(\tilde x,\tilde y)}$ and the
coordinate $\mathbf{x} \in \R$.

If $\mathbf{x}, \tilde x > 0$, we are on the one-dimensional part and
we have a unique geodesics from $(\tilde x,\tilde y)$ to $(\mathbf{x},y)$ and 
along such geodesics we have
\[
 \frac{\d\mu_t}{\d\mm} = \frac{1}{t}\frac{\d\mu_1}{\d\mm} \le \frac{1}{t^3}\frac{\d\mu_1}{\d\mm}\qquad \text{for all }t \in (0,1].
\]

\begin{figure}
  \psfrag{x}{$(\mathbf{x},y)$}
  \psfrag{y}{$(\tilde x, \tilde y)$}
  \psfrag{z}{$\frac{\tilde x}2$}
  \centering
   \includegraphics[width=0.9\textwidth]{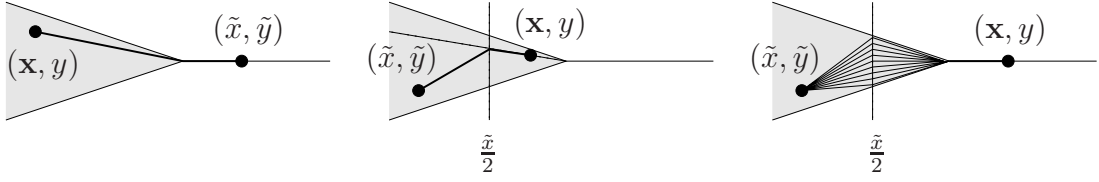}
  \caption{The cases where the selection of geodesics is not euclidean.}
  \label{fig:selection1}
 \end{figure}

If $\mathbf{x} < 0 < \tilde x$, we can transport using a set of geodesics such that each image of a geodesic
is a subset of a set of the form
\[
 \{(x,y)\,:\, y = ax, x \le 0\} \cup L \qquad \text{where }  |a| \le \frac23,
\]
in other words, we first transport along a euclidean geodesic to the origin and then continue along $L$.
See the first case in Figure \ref{fig:selection1}.
Then the induced optimal transport satisfies
\[
 \frac{\d\mu_t}{\d\mm} \le \frac{1}{t}\frac{\d\mu_1}{\d\mm} \le \frac{1}{t^3}\frac{\d\mu_1}{\d\mm}\qquad \text{for all }t \in (0,1].
\]
Compared to the previous case the first equality has now change to an inequality due to the fact that the density might drop when we
pass through the origin.

Next we have to check the different cases when $\tilde x < 0$. Let us first look at the easier case where
$\mathbf{x} < \tilde x$. In this case we can again
use euclidean geodesics. Since the density of the reference measure $\mm$ with respect to $\mathcal{L}_2$ is increasing in $C$
as $x$ increases, we have the estimate
\[
 \frac{\d\mu_t}{\d\mm} \le \frac{1}{t^2}\frac{\d\mu_1}{\d\mm} \le \frac{1}{t^3}\frac{\d\mu_1}{\d\mm}\qquad \text{for all }t \in (0,1].
\]

The final case is when $\tilde x < \mathbf{x}$. In this case we still have three subcases to select 
the geodesic(s) from $(\tilde x, \tilde y)$ to a point $(\mathbf{x},y)$.
The division to subcases at $\frac{\tilde x}{2}$ is due to the fact that all geodesics from  $(\tilde x, \tilde y)$ 
to $(\{\frac{\tilde x}{2}\}\times \R) \cap C$ have the same length.

In the subcase $\tilde x < \mathbf{x}\le \tilde x / 2$ we simply take the euclidean geodesic.
In the subcase $\tilde x / 2 < \mathbf{x} < 0$ we take a geodesic that first goes along the euclidean geodesic 
from $(\tilde x, \tilde y)$ to $(\frac{\tilde x}{2}, \frac{\tilde x y}{2 \mathbf{x}})$ and then along the euclidean geodesic from 
$(\frac{\tilde x}{2}, \frac{\tilde x y}{2 \mathbf{x}})$ to $(\mathbf{x},y)$. See the second case in Figure \ref{fig:selection1}.
The third subcase is when $\mathbf{x} >0$. In this subcase we first 
spread the measure uniformly from $(\tilde x, \tilde y)$ to $\{\tilde x / 2\} \times [-\tilde x /6, -\tilde x /6]$, then contract
along the euclidean geodesics to the origin and then move to $(\mathbf{x},y)$ along the euclidean geodesic.
See the third case in Figure \ref{fig:selection1}.
Now in the case $\tilde x < \mathbf{x} \le \tilde x /2$ we have
\[
 \frac{\d\mu_t}{\d\mm} = \frac{1+\frac{2(\tilde x-\mathbf{x})}{\tilde x}(1-t)}{t^2}\frac{\d\mu_1}{\d\mm} \le \frac{2-t}{t^2}\frac{\d\mu_1}{\d\mm} \le \frac{1}{t^3}\frac{\d\mu_1}{\d\mm}\qquad \text{for all }t \in (0,1].
\]
In the case $\tilde x /2 < \mathbf{x}$ we have
\[
 \frac{\d\mu_t}{\d\mm} = \frac1{t}\frac{\d\mu_1}{\d\mm} \le \frac{1}{t^3}\frac{\d\mu_1}{\d\mm}\qquad \text{for all }t \in [\tilde x /(2\tilde x - 2\mathbf{x}),1]
\]
and
\[
 \frac{\d\mu_t}{\d\mm} = \frac{2-\frac{2(\tilde x-\mathbf{x})}{\tilde x}t}{\left(\frac{2(\tilde x-\mathbf{x})}{\tilde x}t\right)^2}\frac{\d\mu_1}{\d\mm} \le \frac{2-t}{t^2}\frac{\d\mu_1}{\d\mm} \le \frac{1}{t^3}\frac{\d\mu_1}{\d\mm}\qquad \text{for all }t \in (0,\tilde x /(2\tilde x - 2\mathbf{x})].
\]
All cases now being checked we have shown that $(X, \sfd, \mm)$ satisfies $\MCP(0,3)$.

\begin{remark}
 One can check that the space also satisfies the convexity inequality of the $\CD(0,N)$ condition for many transports.
 This is clearly the case for measures transported inside $L$. For transports between $C$ and $L$ this is less trivial,
 but still true. However, for transports in $C$ the $\CD(0,N)$ condition (in fact $\CD(K,\infty)$ for any $K \in \R$) fails.
 The problematic transports are the ones going in the diagonal direction where one is forced to use euclidean geodesics.
 Here the density change of the reference measure $\mm$ with respect to $\mathcal{L}_2$ destroys the $\CD(0,N)$ property.

 It is not clear if the above example satisfies the $\MCP(0,3)$ condition as defined by Sturm \cite{S2006II}. In Sturm's definition
 one should be able to find a distribution of geodesics between any two points so that the $\MCP(K,N)$ condition is satisfied
 along the geodesics given by the distributions. The key point in the definition is that the distribution is symmetric.
\end{remark}

\subsection{Failure of topological maximal diameter theorem}

\begin{figure}
  \psfrag{0}{$0$}
  \psfrag{a}{$\frac14$}
  \psfrag{b}{$-\frac14$}
  \psfrag{c}{$\frac\pi{2}$}
  \psfrag{d}{$-\frac{\pi}2$}
  \centering
   \includegraphics[width=0.6\textwidth]{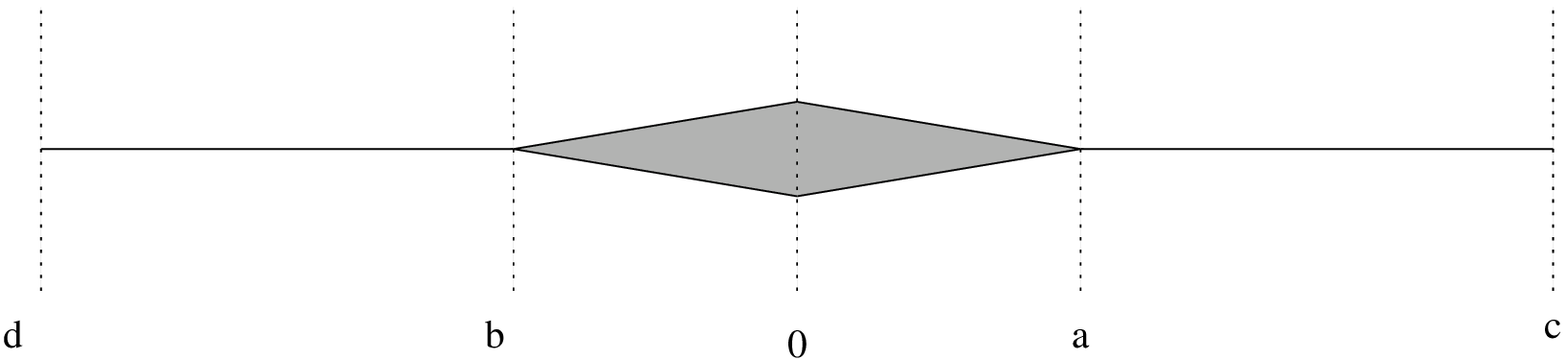}
  \caption{The space $(Y,\sfd,\nn)$ satisfying $\MCP(2,3)$.}
  \label{fig:example2}
 \end{figure}

The space $(Y,\sfd,\nn)$ we describe here with the maximal diameter $\pi$ is similar to our previous example. As a subset of $\R^2$ the space $Y$
is defined as
\[
 Y = D \cup L,
\]
where
\[
 D = \{(x,y) \in \R^2 \,:\, 9|y| \le 1/4 - |x|\} \qquad\text{and}\qquad L = \left([-\pi/{2},-1/4] \cup [1/4,\pi/{2}]\right) \times \{0\}.
\]
See Figure \ref{fig:example2} for an illustration.
Again the distance $\sfd$ is the $l^\infty$-distance \eqref{eq:dist}.
The measure is similar to the measure $\mm$, but this time the projection of the measure to the interval $[-\pi/2,\pi/2]$
is of the form $\d\nn = \cos^2(x)\d\mathcal{L}_1$. Again we define the measure in such a way that the density is constant on $(\{x\} \times \R) \cap Y$
for all $x \in [-\pi/2,\pi/2]$.
The projected space $(\p_1Y, \sfd_e, (\p_1)_\sharp\nn)$ clearly satisfies $\CD(2,3)$.
It is not that surprising that $(Y,\sfd,\nn)$ satisfies $\MCP(2,3)$: the projection deals with the long distances and on
the other hand it is already known that sufficiently small balls in $\R^n$ with the $l^\infty$ norm satisfy $\MCP(n,n+1)$, see \cite{S2006II}.
The transition from the one-dimensional part to the two-dimensional part can then be carried out similarly as in the previous example.

Nevertheless, let us check that the space $(Y,\sfd,\nn)$ really satisfies $\MCP(2,3)$. Let $A \subset Y$ with $\nn(A)>0$
and $(\tilde x, \tilde y) \in Y$. We again transport linearly in the horizontal direction and so it suffices
to check that for all $\mathbf{x} \in [-\pi/2,\pi/2]$ we can select geodesics $\gamma$ from $(\tilde x, \tilde y)$
to $(\{\mathbf{x}\} \times \R) \cap A$ such that
\begin{equation}\label{eq:ex2goal}
 \frac{\d\mu_t}{\d\nn}(\gamma_t) \le \frac{\sin^2(l(\gamma))}{t\sin^2(tl(\gamma))} \frac{\d\mu_1}{\d\nn}(\gamma_1) \qquad \text{for all }t \in (0,1]. 
\end{equation}

We again consider different cases separately depending on the position of $(\tilde x, \tilde y)$ and $\mathbf{x}$.
First of all, if $1/4 \le |\mathbf{x}|, |\tilde x| \le \pi/2$
we can take the measures $\mu_t$ to be uniformly distributed on the slices $(\{x\} \times \R) \cap Y$.
Since the space $(\p_1Y, \sfd_e, (\p_1)_\sharp\nn)$ satisfies $\CD(2,3)$, we have the desired inequality \eqref{eq:ex2goal}
in this case.

The second case is when $|\mathbf{x}| \le 1/4$ and $1/4 \le |\tilde x| \le \pi/2$. In this case we
take the geodesics along which we transport to be such that their images are subsets of sets of the type 
\[
L \cup \left\{(x,y) \in \R^2 \,:\, cy = 1/4 - |x| \ge 0\right\}, \qquad \text{with } |c| \ge 9.
\]
Because inside $D$ we are changing the density as if we were transporting in $(\p_1Y, \sfd_e, (\p_1)_\sharp\nn)$
and the density relative to this transport only drops when we leave $D$, again the fact that we have the inequality 
\eqref{eq:ex2goal} follows from the observation that the space $(\p_1Y, \sfd_e, (\p_1)_\sharp\nn)$ satisfies $\CD(2,3)$.

The final case to check is when $|\tilde x| \le 1/4$. By symmetry we may assume that $\tilde x \le \mathbf{x}$.
As in the previous example we define the chosen geodesics from $(\tilde x, \tilde y)$ to a point $(\mathbf{x},y) \in A$ in three subcases.
This time the division to subcases is driven by the fact that the distances from $(\tilde x, \tilde y)$
to all the points in $(\{(1/4+4\tilde x)/5\}\times\R) \cap D$ are equal, meaning that beyond that distance we can
transport along vertically uniformly distributed measures.

Let us list the subcases.
First, if $\tilde x < \mathbf{x} \le (1/4+4\tilde x)/5$, we take the euclidean geodesics.
In the subcase $(1/4+4\tilde x)/5 < \mathbf{x} < 1/4$ we take a geodesic that first goes along the euclidean geodesic 
that goes via the point $(1/4,0)$
from $(\mathbf{x},y)$ to $(\max\{0,(1/4+ 4\tilde x)/5\}, y')$ and then along the euclidean geodesic from 
$(\max\{0,(1/4+ 4\tilde x)/5\}, y')$ to $(\tilde x, \tilde y)$.
The third subcase is when $\mathbf{x} >1/4$. In this subcase we first transport the mass along the euclidean
geodesic to the point $(1/4,0)$, then spread the measure uniformly to $(\{\max\{0,(1/4+ 4\tilde x)/5\}\} \times \R) \cap D$ along
euclidean geodesics and finally contract it to $(\tilde x, \tilde y)$ along euclidean geodesics.

Let us check the density in the different subcases. To deal with the first subcase $\tilde x < \mathbf{x} \le (1/4+4\tilde x)/5$
we observe that $h_t \le (5/4-t/4)h_1$ for the vertical height $h_t$ of $D$ at $\gamma_t$, where $\gamma$ is a geodesic from $(\tilde x, \tilde y)$
to $(\mathbf{x},y)$.
From this (and the estimates below) we have that
\begin{equation}\label{eq:essential}
 \frac{\d\mu_t}{\d\nn} \le \frac{h_t}{t^2h_1}\frac{\d\mu_1}{\d\nn} \le \frac{(5/4-t/4)}{t^2}\frac{\d\mu_1}{\d\nn} 
  \le \frac{\sin^2(l)}{t\sin^2(tl)} \frac{\d\mu_1}{\d\nn}
\end{equation}
for all $t \in (0,1]$ and  $l \in (0,1/2]$ and so \eqref{eq:ex2goal} follows.
(To see the third inequality in \eqref{eq:essential} define $f_l(t) := \left(\frac54-\frac{t}4\right)\sin^2(tl)-t\sin^2(l)$ and notice
 that $f_l(t) \le 0$ for all $t \in [0,1]$ and $l \in (0,1/2]$ since $f_l(0)=f_l(1)=0$
 and $f_l''(t) = -\frac{l}{2}\sin(2tl) + (5-t)\frac{l^2}{2}\cos(2tl) \ge 0$ for all $t \in [0,1]$ and $l \in (0,1/2]$.)
 
In the second and third subcases $(1/4+4\tilde x)/5 < \mathbf{x}$, until we hit $(\max\{0,(1/4+ 4\tilde x)/5\}, y')$ 
we have \eqref{eq:ex2goal} by the projection as before.
After $(\max\{0,(1/4+ 4\tilde x)/5\}, y')$ we have \eqref{eq:ex2goal} by the estimate \eqref{eq:essential}
for $l \in (0,\frac12]$.
When $l > \frac12$ we can estimate the time $t_0$ when we hit $(\max\{0,(1/4+ 4\tilde x)/5\}, y')$
by
\[
 t_0 \le \frac1{l}\left(\frac{\frac14+4\tilde x}{5} - \tilde x\right) \le 2\frac{\frac14-\tilde x}{5} \le \frac15.
\]
In this case obtaining \eqref{eq:ex2goal} is easy, since verifying \eqref{eq:essential}
for $t \in [0,\frac15]$ and $l \in (\frac12,\frac{\pi}{2}+\frac14)$ reduces to verifying that $f_l(t) \le 0$ for such $t$ and $l$. 
This follows for example by estimating
\[
 \left(\frac54-\frac{t}4\right)\sin^2(tl)\le \frac54 \sin^2(tl) \le \frac54(tl)^2 \le \frac14 tl^2\le t\sin^2(l)
\]
for $t \in [0,\frac15]$ and $l \in (\frac12,\frac{\pi}{2}+\frac14)$.
This completes the proof of the $\MCP(2,3)$ condition for our space $(Y,\sfd, \nn)$.
     
\subsection{Final remarks}\label{sec:remarks}

Both of our examples above had boundary points in the sense that there are geodesics
that cannot be infinitely extended as local geodesics. It is natural to ask if one can modify the examples to
have no boundary. By identifying the top and bottom boundaries in the two-dimensional parts of the spaces
we should be able to get rid of most of the boundary points. However, even after the identification
the point where the space changes from one-dimensional to two-dimensional will be a boundary point in the above sense.
The nonextendable geodesics are the nonvertical geodesics in the two-dimensional side that are travelling more in the
vertical direction than in the horizontal direction.
When they reach the point where the dimension changes they cannot be extended to the one-dimensional part so that they
are geodesic in a neighbourhood of the point where the dimension changes.

It is also interesting to consider the tangent spaces of the spaces $(X,\sfd,\mm)$ and $(Y,\sfd,\nn)$
in the sense of measured Gromov-Hausdorff limits of blow-ups. First observation is that the space
$(X,\sfd,\mm)$ is the mGH-tangent of $(Y,\sfd,\nn)$ at the points $(a,0)$ and $(-a,0)$. This way
we could prove the fact that $(X,\sfd,\mm)$ satisfies $\MCP(0,3)$ by proving that any tangent of
a $\MCP(K,3)$ space does.

Notice that the tangents of $(X,\sfd,\mm)$ at points $X \setminus \{(0,0)\}$ and $(Y,\sfd,\nn)$
at points $Y\setminus \{(1/4,0),(-1/4,0)\}$ are all convex subsets of normed spaces.
On the other hand we know from \cite{J2009} that some subRiemannian spaces like the Heisenberg group satisfy $\MCP(K,N)$.
With only such nice examples at hand it is natural to ask what are the tangent spaces of
$\MCP(K,N)$ spaces (almost everywhere) in general? This is also related to the uniqueness of the tangents,
compare to the following result by Le Donne \cite{LD2011}.

\begin{theorem}[Le Donne]
 Let $(Z, \sfd, \mu)$ be a doubling-measured geodesic metric space. Assume that for $\mu$-almost every
 $x \in Z$ the space $Z$ has only one tangent. Then for $\mu$-almost every point $x \in Z$ the tangent
 space is a Carnot group $\mathbf{G}$ endowed with a subFinsler left-invariant metric with the first
 layer of the Lie algebra of $\mathbf{G}$ as horizontal distribution.
\end{theorem}

\end{document}